\documentclass[11pt]{amsart}

\newcommand{\private}[1]{}
\usepackage{amssymb, amsmath, amscd, amsthm, color, epsfig,url}

\usepackage[all]{xy}          
\xyoption{dvips}              

\makeatletter
\renewcommand\l@subsection{\@tocline{2}{0pt}{2pc}{5pc}{}}
\makeatother

\newcommand{\hofiber}{\operatorname{hofiber}}

\newcommand{\hocofiber}{\operatorname{hocofiber}}
\newcommand{\holim}{\operatorname{holim}}

\newcommand{\hocolim}{\operatorname{hocolim}}

\newcommand{\Map}{\operatorname{Map}}

\newcommand{\Link}{\operatorname{Link}}

\theoremstyle{plain}
\newtheorem{thm}{Theorem}[section]

\theoremstyle{definition}

\newtheorem{def/ex}[thm]{Definition/Example}

\theoremstyle{remark}

\newcommand{\refS}[1]{Section~\ref{S:#1}}
\newcommand{\refT}[1]{Theorem~\ref{T:#1}}

\begin{document}


\title[A Stable Range Description of the Space of Link Maps]{A Stable Range Description of the Space of Link Maps}


\author{Thomas G. Goodwillie}
\address{Department of Mathematics, Brown University, Providence, RI}
\email{tomg@math.brown.edu}

\author{Brian A. Munson}
\address{Department of Mathematics, Wellesley College, Wellesley, MA}
\email{bmunson@wellesley.edu}
\urladdr{http://palmer.wellesley.edu/\~{}munson}




\begin{abstract}
{\bf Version: \today}
We study the space $\Link(P,Q; N)$ of link maps, maps from $P\coprod Q$ to $N$ such that the images of $P$ and $Q$ are disjoint. We give a range of dimensions, interpreted as the connectivity of a certain map, in which the cobordism class of the ``linking manifold'' is enough to distinguish the homotopy class of one link map from another.
\end{abstract}

\maketitle

\tableofcontents

\parskip=4pt
\parindent=0cm


\section{Introduction}\label{S:Intro}


Let $N$ be a smooth manifold and let $P$ and $Q$ be smooth compact manifolds. A (smooth)  \emph{link map} of $P$ and $Q$ in $N$ is a pair $(f:P\rightarrow N, g:Q\rightarrow N)$ of smooth maps such that $f(P)$ is disjoint from $g(Q)$. The set of link maps, denoted by $\Link(P,Q;N)$, is an open subspace of $\Map(P,N)\times \Map(Q,N)=\Map(P\coprod Q,N)$.

For brevity here we will write $\mathcal M$ for $\Map(P,N)\times \Map(Q,N)$ and denote the complement of the set of link maps in $\mathcal M$ by $\mathcal B$. We prove that a certain ``linking number'' map

\begin{equation*}
\ell:\hofiber_{(f_1,g_1)}(\mathcal{M}-\mathcal B\rightarrow \mathcal{M})\rightarrow \Omega Q^{TN-(TP\oplus TQ)}_+\holim(P\stackrel{f_1}{\rightarrow} N\stackrel{g_1}{\leftarrow} Q)
\end{equation*}

is $(2(n-p-q)-3)$-connected, where $p$, $q$, and $n$ are the dimensions of the manifolds. The map was defined by the second author in \cite{M:LinkNumber}, although the version we reference below is a homotopy-theoretic one give in \cite{KlWi:HIT1}. Its domain is the homotopy fiber of the inclusion $\mathcal M-\mathcal B\rightarrow\mathcal M$ with respect to any point $(f_1,g_1)\in \mathcal M-\mathcal B$. Its codomain is the infinite loopspace associated to the Thom spectrum of a virtual vector bundle. Both of these spaces are $(n-p-q-2)$-connected. In the case when $p+q=n-1$ it was shown in  \cite{M:LinkNumber} that  the effect of the map $\ell$ on $\pi_0$ can be interpreted as a generalized linking number.

Remark: Functor calculus (the manifold version developed in \cite{W:EI1,GW:EI2}) offers one point of view on link maps. Consider the functor $(U,V)\mapsto \Link(U,V;N)$ whose domain is the poset $\mathcal{O}(P\coprod Q)=\mathcal{O}(P)\times\mathcal{O}(Q)$ of open subsets of $P\coprod Q$. Its best linear approximation is $\Map(U,N)\times\Map(V,N)$. Our result can be interpreted as a statement about a quadratic approximation to the same functor, but we will not pursue this here. This work overlaps the recent work of Klein and Williams \cite{KlWi:HIT1}. In particular, some of the material in \refS{codim} also appears in \cite{KlWi:HIT1}.

Our main result is:

\begin{thm}\label{T:linkstable}
The map $$\Lambda:\Sigma\hofiber(\mathcal M-\mathcal B\rightarrow \mathcal M)\rightarrow Q^{TN-(TP\oplus TQ)}_+\holim(P\rightarrow N\leftarrow Q)$$
adjoint to $\ell$ is $(2(n-p-q)-1)$-connected. 
\end{thm}

The fact that $\ell$ is $(2(n-p-q)-3)$-connected then follows immediately by the Freudenthal theorem, since the domain of $\ell$ is $(n-p-q-1)$-connected. Note that the connectivity claimed for $\Lambda$ is negative if $p+q\geq n$, so it is no loss to assume $p+q<n$.

\subsection{Conventions}\label{S:conventions}

A space $X$ is $k$-connected if for every $j$ with $-1\leq j\leq k$ every map $S^j\rightarrow X$ can be extended to a map $D^{j+1}\rightarrow X$. (In other words, $(-1)$-connected means nonempty and if $k\geq 0$ then $k$-connected means that there is exactly one path-component and that the homotopy groups vanish through dimension $k$.)
A map is $k$-connected if each of its homotopy fibers is $(k-1)$-connected. A (weak)  \textsl{equivalence} is an $\infty$-connected map.

We write $QX=\Omega^{\infty}\Sigma^{\infty}X$ if $X$ is a based space. If $X$ is unbased, then $X_+$ means $X$ with a disjoint basepoint added and $Q_+X$ means $Q(X_+)$. For a vector bundle $\xi$ over a space X, the unit disk bundle and the unit sphere bundle are $D(X;\xi)$ and $S(X;\xi)$. The Thom space $X^\xi$ is the quotient $D(X;\xi)/S(X;\xi)$, or equivalently the homotopy cofiber of the projection $S(X;\xi)\rightarrow X$. If $\xi$ and $\eta$ are two vector bundles on $X$, then by choosing a vector bundle monomorphism $\eta\rightarrow \epsilon^i$ to a trivial bundle we can define $Q_+^{\xi-\eta}X=\Omega^iQX^{\xi\oplus\epsilon^i/\eta}$. This is essentially independent of the choice of $i\geq 0$ and vector bundle monomorphism, in the sense that for large $i$ the weak homotopy type of this space is independent of those choices.

\section{Sketch of the proof of \refT{linkstable}}\label{S:proof}

To prove \refT{linkstable} we will use the diagram (\ref{bigdiagram}) below and obtain the connectivity of $\Lambda$ from the connectivities of all the other maps. For this we must introduce another closed set $\mathcal V\subset \mathcal B$. Recall that a point $(f,g)\in \mathcal{M}$ belongs to $\mathcal B$ if  the statement $f(x)=z=g(y)$ holds for some pair $(x,y)\in P\times Q$ and some point $z\in N$. The closed set $\mathcal B$ has codimension $n-p-q$ in $\mathcal M$ in some sense. Inside this space of ``bad'' maps is a set $\mathcal{V}$ of ``very bad'' maps, having codimension in $\mathcal M$ is $2(n-p-q)$. A point $(f,g)$ is in $\mathcal V$ if either the statement $f(x)=z=g(y)$ holds for more than one choice of $(x,y,z)$ or else it holds for one such choice in such a way that the associated map of tangent spaces $T_xP\oplus T_yQ\rightarrow T_zN$ is not injective. The set $\mathcal B-\mathcal V$ may be regarded as a submanifold of $\mathcal{M}$. It has maps to $P$, $Q$, and $N$ given by $x$, $y$, and $z$. Pulling back tangent bundles via these maps, we obtain vector bundles on $\mathcal B-\mathcal V$, which we will denote simply by $TP$, $TQ$, and $TN$. There is also a monomorphism $TP\oplus TQ\rightarrow TN$, and its cokernel $TN/(TP\oplus TQ)$ may be thought of as the normal bundle of $\mathcal B-\mathcal V$ in $\mathcal M$. 


The next result immediately implies \refT{linkstable}.

\begin{thm}\label{T:connstatements}
In the diagram below, the maps $F$ and $H$ are equivalences, the maps $G,C,$ and $D$ are $(2(n-p-q)-1)$-connected, and the map $E$ is $(3(n-p-q)-2)$-connected.
\end{thm}

\begin{equation}\label{bigdiagram}
\xymatrix{
\Sigma\hofiber(\mathcal M-\mathcal B\rightarrow \mathcal{M})\ar[r]^{\Lambda} & Q_+^{TN-(TP\oplus TQ)}\holim(P\rightarrow N\leftarrow Q) \\
 & Q_+^{TN-(TP\oplus TQ)}\hofiber(\mathcal B -\mathcal V\rightarrow \mathcal{M})\ar[u]_D\\
\Sigma\hofiber(\mathcal M-\mathcal B\rightarrow \mathcal M- \mathcal V)\ar[uu]^G & Q\hofiber(\mathcal B-\mathcal V\rightarrow \mathcal{M})^{TN/(TP\oplus TQ)}\ar[u]_H\\
\hofiber(\mathcal B-\mathcal V\rightarrow \mathcal M-\mathcal V)^{TN/(TP\oplus TQ)} \ar[r]^E\ar[u]^F & \hofiber(\mathcal B-\mathcal V\rightarrow \mathcal{M})^{TN/(TP\oplus TQ)}\ar[u]_C\\
}
\end{equation}

We now briefly define the maps in the diagram and explain about their connectivities. Steps that are sketchy here will be filled in in the following sections. Let $c=n-p-q$.

The equivalence $F$ is essentially an instance of the following general fact. If $Y$ is a smooth submanifold of $X$ and also a closed subset, then the suspension of the homotopy fiber of the inclusion $X-Y\rightarrow X$ is equivalent to the Thom space, over the homotopy fiber of $Y\rightarrow X$, of the normal bundle of $Y$ in $X$. This general fact will be proved, and adapted to the present function-space setting, in \refS{bundles&cofibers}.

$G$ is an inclusion map. Since $\mathcal V$ has codimension $2c$ in $\mathcal M$, the inclusion $\mathcal M-\mathcal V\rightarrow \mathcal M$ is $(2c-1)$-connected. (This will be worked out in detail in \refS{codim}) Therefore the map of homotopy fibers is $(2c-2)$-connected and the map $G$ of suspensions is $(2c-1)$-connected. 

$E$ is a map of Thom spaces. For a $k$-connected map $Z\rightarrow W$ of spaces and a vector bundle $\xi$ on $W$ with fiber dimension $d$, the associated map $Z^\xi\rightarrow W^\xi$ is $(k+d)$-connected. In our case $d=c$ and $k=2c-2$; the inclusion of $\hofiber(\mathcal B-\mathcal V\rightarrow \mathcal M-\mathcal V)$ into $\hofiber(\mathcal B-\mathcal V\rightarrow \mathcal{M})$ is $(2c-2)$-connected, again because the inclusion of $\mathcal M-\mathcal V$ into $\mathcal M$ is $(2c-1)$-connected.
 
$C$ is the canonical map $Z\rightarrow QZ$, where the space $Z$ is $(c-1)$-connected, being the Thom space of a vector bundle of rank $c$. By the Freudenthal Theorem, the map is $(2c-1)$-connected.

The equivalence $H$ is simply a matter of rewriting the Thom spectrum of a virtual vector bundle $\xi -\eta$ as the suspension spectrum of the Thom space of $\xi/\eta$ when $\eta$ is a subbundle of $\xi$. 

The map $D$ arises from a $(c-1)$-connected map from $\hofiber(\mathcal B-\mathcal V\rightarrow \mathcal{M})$ to $\holim(P\rightarrow N\leftarrow Q)$. To explain further, we introduce the space $\mathcal{\tilde B}$ of all 
$((f,g),x,y,z)\in \mathcal M\times P\times Q\times N$ 
such that $f(x)=z=g(y)$. Projection to $\mathcal M$ gives a map from $\mathcal {\tilde B}$ onto $\mathcal B$. Let $\mathcal {\tilde V}\subset \mathcal{\tilde B}$ be the preimage of $\mathcal V$. The projection $\mathcal {\tilde B}-\mathcal {\tilde V}\rightarrow \mathcal B-\mathcal V$ is an isomorphism. The inclusion $\mathcal {\tilde B}-\mathcal {\tilde V}\rightarrow \mathcal {\tilde B}$ is $(c-1)$-connected for reasons of codimension (again, the details are in section 3), and therefore the induced map $\hofiber(\mathcal B-\mathcal V\rightarrow \mathcal{M})\rightarrow \hofiber(\mathcal {\tilde B}\rightarrow \mathcal{M})\simeq\holim(P\rightarrow N\leftarrow Q)$ is also $(c-1)$-connected. There are vector bundles $TP$, $TQ$, and $TN$ on $ \mathcal {\tilde B}$ pulling back to their namesakes on $\mathcal B-\mathcal V$. (The monomorphism $df\oplus dg:TP\oplus TQ\rightarrow TN$ is not available on the $\holim(P\rightarrow N\leftarrow Q)$ side,  which is why we switched from Thom spaces to Thom spectra).

We end this section with a brief account of the commutativity of diagram (\ref{bigdiagram}). First we need to define the map $\Lambda$. As mentioned in \refS{Intro}, $\Lambda$ is adjoint to a map $\ell:\hofiber_{(f_1,g_1)}(\mathcal M-\mathcal B\to\mathcal M)$, which is a composite defined as follows. Let $(f_t,g_t)\in\hofiber_{(f_1,g_1)}(\mathcal M-\mathcal B\to\mathcal M)$. The map $\mathcal M\to\Map(P\times Q,N\times N)$ given by $(f,g)\mapsto f\times g$ induces a map $\hofiber_{(f_1,g_1)}(\mathcal M-\mathcal B\to\mathcal M)\to\hofiber_{f_1\times g_1}(\Map(P\times Q,N\times N-\Delta_N)\to \Map(P\times Q,N\times N))$. We can identify the latter homotopy fiber as a space of sections as follows.

Let $E=\holim(P\times Q\stackrel{f_1\times g_1}{\longrightarrow}N\times N\longleftarrow N\times N-\Delta_N)$. The projection map $E\to P\times Q$ is a fibration with fiber over $(p,q)$ the space $\Phi_2(N)=\hofiber_{(f_1(p),g_1(q)}(N\times N-\Delta_N\to N\times N)$. Let $\Gamma(P\times Q,E)$ its space of sections. This space of sections has a preferred basepoint given by $(f_1,g_1)$. It is equivalent to $\hofiber_{f_1\times g_1}(\Map(P\times Q,N\times N-\Delta_N)\to\Map(P\times Q, N\times N))$ by inspection. Let $Q_fS_fE\to P\times Q$ be the fibration whose fibers are $QS\Phi_2(N)$, where $S$ stands for the unreduced suspension. The canonical map $\Gamma(P\times Q,E)\to\Omega\Gamma(P\times Q,Q_fS_fE)$ is easily shown to be $(2n-p-q-1)$-connected, and there is an equivalence $\Omega\Gamma(P\times Q,Q_fS_fE)\simeq\Omega Q_+^{TN-(TP\oplus TQ)}\holim(P\times Q\stackrel{f_1\times g_1}{\to}N\times N\leftarrow\Delta_N)$ which is the identity on the loop coordinate. Moreover, there is a homeomorphism $\holim(P\times Q\stackrel{f_1\times g_1}{\to}N\times N\leftarrow\Delta_N)\cong \holim(P\stackrel{f_1}{\longrightarrow}N\stackrel{g_1}{\longleftarrow}Q)$. The composite map $\hofiber_{(f_1,g_1)}(\mathcal M-\mathcal B\to\mathcal M)\to\Omega Q_+^{TN-(TP\oplus TQ)}\holim(P\stackrel{f_1}{\longrightarrow}N\stackrel{g_1}{\longleftarrow}Q)$ is the map $\ell$, and $\Lambda$ is its adjoint. See \cite[Section 9]{KlWi:HIT1}.

Now let $(f_t,g_t,v)\in\hofiber(\mathcal B-\mathcal V\to\mathcal M-\mathcal V)^{TN/TP\oplus TQ}$. Here $v$ is a vector of length $0\leq|v|\leq1$, and $(f_t,g_t,v)$ is identified to a point when $|v|=1$. After applying the maps $E,C,H$, and $D$ in diagram (\ref{bigdiagram}), it is clear that $(f_t,g_t,v)$ is sent to $((x,\beta, y),v)\in Q_+^{TN-(TP\oplus TQ)}\holim(P\to N\leftarrow Q)$, where $(x_0,y_0)\in P\times Q$ is the unique pair such that $f_0(x_0)=g_0(y_0)$ and $\beta:I\to N$ is the path defined by $\beta(s)=f_{1-2s}(x_0)$ for $0\leq s\leq1/2$ and $\beta(s)=g_{2s-1}(y_0)$ for $1/2\leq s\leq 1$.

Now we must apply $F,G$, and $\Lambda$ to $(f_t,g_t,v)$. A careful examination of the material in \refS{bundles&cofibers} reveals that $F$ sends $(f_t,g_t,v)$ to the point $s\wedge (\tilde{f}_t,\tilde{g}_t)$, where $s=1-|v|$ and $(\tilde{f}_t,\tilde{g}_t)\in\hofiber(\mathcal M-\mathcal B\to\mathcal M)$ is defined as follows. For $s\leq t\leq 1$, $(\tilde{f}_t,\tilde{g}_t)=(f_{(t-s)/(1-s)},g_{(t-s)/(1-s)})$. For $0\leq t\leq s$, $(\tilde{f}_t,\tilde{g}_t)$ has the following properties: $(\tilde{f}_t,\tilde{g}_t)\in \mathcal M-\mathcal B$ for $t<s$, $(\tilde{f}_s,\tilde{g}_s)=(f_0,g_0)$ has a unique pair $(x_0,y_0)\in P\times Q$ such that $f_0(x_0)=g_0(y_0)=z_0\in N$ and such that $f'_0(x_0)-g'_0(y_y)\in T_{z_0}N$, when projected to $T_{z_0}/T_{x_0}P\oplus T_{y_0}Q$, is equal to $v$ (here $f'_0$ and $g'_0$ are the derivatives with respect to $t$). From this description of $F$ and the description of $\Lambda$ above, the diagram commutes.


\section{Codimension and Connectivity}\label{S:codim}


The proof outlined above uses that the pair $(\mathcal M,\mathcal M-\mathcal V)$ is $(2n-2p-2q-1)$-connected and that the pair
$(\mathcal {\tilde B},\mathcal {\tilde B}-\mathcal {\tilde V})$ is $(n-p-q-1)$-connected. We now justify these statements more carefully.

For the first, it suffices if for every smooth manifold $W$ of dimension $j<2n-2p-2q$, for every map of pairs $\phi:(W,\partial W)\rightarrow (\mathcal M,\mathcal M-\mathcal V)$, there is a homotopy of pairs to a map that is disjoint from $\mathcal V$. 

Consider the adjoint map $\Phi:W\times (P\coprod Q)\rightarrow N$. By a preliminary homotopy we can assume that $\Phi$ is smooth, and we can make the homotopy small enough in the $C^0$ sense so that it corresponds to a homotopy of pairs. If we can show that the condition $\phi^{-1}(\mathcal V)=\emptyset$ holds for a dense set of all such smooth maps $\Phi$, then another small homotopy will complete the job. For the density statement we will use the multijet transversality theorem of Mather \cite[Proposition 3.3] {Ma:StabilityV} (which appears in \cite{GG:stablemaps} as Theorem 4.13).

Recall the setup: Two smooth maps $\Phi,\Psi:X\rightarrow Y$ have the same $m$-jet at $x\in X$ if $\Phi(x)=\Psi(x)$ and $\Phi$ and $\Psi$ have the same derivatives through order $m$. Let $X^{(r)}\subset X^r$ be the space of configurations of $r$ distinct points in $X$. The maps $\Phi$ and $\Psi$ have the same $m$-multijet at $(x_1,\dots ,x_r)\in X^{(r)}$ if for every $i\in \lbrace 1,\dots r\rbrace$ they have the same $m$-jet at $x_i$. The manifold $J^{(r)}_m(X,Y)$ of multijets has a point for each $r$-tuple $(x_1,\dots ,x_r)$ and each equivalence class of maps as above. A smooth map $\Phi:X\rightarrow Y$ determines a smooth map $j^{(r)}_m(\Phi):X^{(r)}\rightarrow J^{(r)}_m(X,Y)$. The multijet transversality theorem asserts that, for every submanifold $Z$ of $J^{(r)}_m(X,Y)$, the set of all $\Phi$ such that $j^{(r)}_m(\Phi)$ is transverse to $Z$ is a countable intersection of dense open sets in the function space $\Map(X,Y)$. It follows that such a set, or even the intersection of countably many such sets, is dense.

We now introduce various submanifolds $Z$ of $J^{(r)}_m(W\times (P\coprod Q),N)$, for various values of $r$ and $m$. The point is that the condition $\phi^{-1}(\mathcal V)=\emptyset$ will hold if and only if for each of these the set $j^{(r)}_m(\Phi)$ is disjoint from $Z$. The codimension of $Z$ will always be big enough so that  in order for $j^{(r)}_m(\Phi)$ to be transverse to $Z$ it must be disjoint. Therefore the theorem will guarantee that there are maps $W\rightarrow \mathcal M-\mathcal V$ arbitrarily close to a given map $\Phi:W\rightarrow \mathcal M$.

Let $k$ be the dimension of $W$. 
We consider the various ways in which $\phi$ could hit $\mathcal V$. 

\begin{enumerate}
\item There might exist distinct $x_1$ and $x_2$ in $P$ and distinct $y_1$ and $y_2$ in $Q$ such that for some $w\in W$ we have $\Phi(w,x_1)=\Phi(w,y_1)$ and $\Phi(w,x_2)=\Phi(w,y_2)$. Then the point 
$$((w,x_1),(w,x_2),(w,y_1),(w,y_2))\in (W\times (P\coprod Q))^{(4)}$$
maps into a certain submanifold of $J^{(4)}_0(W\times (P\coprod Q),N)$ whose codimension is $3k+2n$. (That is $3j$ to make four points of $W$ equal to each other and $2n$ for two coincidences in $N$.) This codimension is greater than the dimension $4k+2p+2q$ of (the relevant open and closed part of) $(W\times (P\coprod Q))^{(4)}$, so that transverse means disjoint.

\item There might exist distinct $x_1$ and $x_2$ in $P$ and $y$ in $Q$ such that $\Phi(w,x_1)=\Phi(w,y)=\Phi(w,x_2)$. This leads to a submanifold of $J^{(3)}_0(W\times (P\coprod Q),N)$ whose codimension is $2k+2n$, greater than the dimension $3k+2p+q$ of (part of) $(W\times (P\coprod Q))^{(3)}$.

\item There might exist $x$ in $P$ and distinct $y_1$ and $y_2$ in $Q$ such that $\Phi(w,x)=\Phi(w,y_1)=\Phi(w,y_2)$. The relevant manifold has codimension $2k+2n$ in $J^{(3)}_0(W\times (P\coprod Q),N)$, greater than $3k+p+2q$.

\item There might exist $x\in P$ and $y\in Q$ such that $\Phi(w,x)=z=\Phi(w,y)$ and such that the linear map $T_xP\oplus T_yQ\rightarrow T_zN$ given by differentiation of $\phi(w)$ at $x$ and $y$ has rank less than $p+q$. For each fixed rank $r<p+q$ this leads to a submanifold of $J^{(2)}_1(W\times (P\coprod Q),N)$ whose codimension $k+(n-r)(p+q-r)$ is greater than $2k+p+q$.
\end{enumerate}

This completes the proof that the pair $(\mathcal M,\mathcal M-\mathcal V)$ is $(2n-2p-2q-1)$-connected.

To prove that the pair
$(\mathcal {\tilde B},\mathcal {\tilde B}-\mathcal {\tilde V})$ is $(n-p-q-1)$-connected, essentially the same kind of standard dimension-counting will succeed, but a simple reference as before to the multijet transversality theorem will not suffice because $\mathcal {\tilde B}$ is not simply the space of maps from one manifold to another. 

First observe that both the projection $\mathcal{\tilde B}\rightarrow P\times Q\times N$
and its restriction $\mathcal {\tilde B}-\mathcal {\tilde V}\rightarrow P\times Q\times N$ are fibrations. It therefore suffices if, for a point $(x_0,y_0,z_0)\in P\times Q\times N$, the pair $(\mathcal {\tilde B}_0,\mathcal {\tilde B}_0-\mathcal {\tilde V}_0)$ of fibers is $(n-p-q-1)$-connected. Here $\mathcal {\tilde B}_0\subset\mathcal M$ is the set of all $\phi$ such that $\phi(x_0)=z_0=\phi(y_0)$, and $\mathcal {\tilde V}_0\subset\mathcal {\tilde B}_0$ is the set of all $\phi$ such that in addition at least one of the following is true:

\begin{enumerate}
\item $\phi(x)=\phi(y)$ for some $x\in P-x_0$ and some $y\in Q-y_0$

\item $\phi(x)=z_0$ for some $x\in P-x_0$

\item $\phi(y)=z_0$ for some $y\in Q-y_0$

\item The linear map $T_{x_0}P\oplus T_{y_0}Q\rightarrow T_{z_0}N$ has rank less than $p+q$. 
\end{enumerate}

To deal first with (4), note that $\mathcal {\tilde B}_0$ is fibered over the space $\mathcal L$ of all linear maps $T_{x_0}P\oplus T_{y_0}Q\rightarrow T_{z_0}N$. Let $\mathcal L^{max}\subset \mathcal L$ be the open set of maps of rank $p+q$ and let $\mathcal {\tilde B}_0^{max}\subset \mathcal {\tilde B}_0$ be its preimage. The pair $(\mathcal {\tilde B}_0, \mathcal {\tilde B}_0^{max})$ is $(n-p-q-2)$-connected (one better than needed), because the pair $(\mathcal L,\mathcal L^{max})$ is $(n-p-q-2)$-connected, because the closed set $\mathcal L-\mathcal L^{max}$ is the union of finitely many submanifolds having codimension at least $n-p-q-1$.

It remains to show that the pair $(\mathcal {\tilde B}_0^{max},\mathcal {\tilde B}-\mathcal {\tilde V})$ is $(n-p-q-1)$-connected. Both $\mathcal {\tilde B}_0^{max}$ and $\mathcal {\tilde B}-\mathcal {\tilde V}$ fiber over $L^{max}$, so we can replace the two spaces by their fibers, say $\mathcal {\tilde B}_L$ and $\mathcal {\tilde B}_L-\mathcal {\tilde V}_L$ over a given $L\in \mathcal L$.

Now given a map $\phi: W\rightarrow \mathcal {\tilde B}_L$, we want to perturb it slightly so as to eliminate behaviors (1), (2) and (3). None of these can occur for $x$ near $x_0$ or $y$ near $y_0$ anyway, given the choice of $L$, so we look for perturbations that are fixed near $x_0$ and $y_0$. In other words, we look for a small compactly supported change in the map $\Phi: W\times ((P-x_0)\coprod (Q-y_0))\rightarrow N$. This goes as before: case (1) leads to  a submanifold of $J^{(2)}_0(W\times ((P-x_0)\coprod (Q-y_0)),N)$ with codimension $j+n$, greater than $2j+p+q$; case (2) leads to  a submanifold of $J^{(1)}_0(W\times ((P-x_0)\coprod (Q-y_0)),N)$ with codimension $n$, greater than $j+p$.; and case (3) leads to  a submanifold of $J^{(1)}_0(W\times ((P-x_0)\coprod (Q-y_0)),N)$ with codimension $n$, greater than $j+q$.

\section{Normal bundles and homotopy cofibers}\label{S:bundles&cofibers}

Suppose that $X$ is a smooth manifold, and that the closed subset $Y\subset X$ is a smooth submanifold with normal bundle $\nu$. 

Of course, the Thom space $Y^{\nu}$ is equivalent to the homotopy cofiber of the inclusion map $X-Y\rightarrow X$. This follows from the fact that there is a homotopy pushout square

\begin{equation}\label{manifoldpushout}
\xymatrix{
S(Y;\nu)\ar[r]\ar[d] & D(Y;\nu)\ar[d]\\
X-Y\ar[r] & X\,.\\
}
\end{equation}

The homotopy fibers over $X$ of the four spaces above form another homotopy pushout square
$$\xymatrix{
\hofiber(S(Y;\nu)\rightarrow X)\ar[r]\ar[d] & \hofiber(D(Y;\nu)\rightarrow X)\ar[d]\\
\hofiber(X-Y\rightarrow X)\ar[r] & \hofiber(X\rightarrow X)\simeq *\,.\\
}
$$
Comparing homotopy cofibers of the rows in this square, we obtain an equivalence
$$
\hofiber(Y\rightarrow X)^\nu\rightarrow \Sigma \hofiber(X-Y\rightarrow X)
$$
Here we have written $\nu$ for the pullback of $\nu$ to $\hofiber(Y\rightarrow X)$. 

We need statements like those above in which the manifolds $X$ and $Y$ are replaced by the function spaces $\mathcal M-\mathcal V$ and $\mathcal B-\mathcal V$ and the role of the normal bundle is played by the vector bundle $TN/(TP\oplus TQ)$ on $\mathcal B-\mathcal V$. The only little difficulty is that the square (\ref{manifoldpushout}) depended on having a tubular neighborhood. We will write down a substitute for (\ref{manifoldpushout}) that avoids this dependence.

Let $P(Y,X)$ be the space of all smooth paths $\gamma:[0,1]\rightarrow X$ such that $\gamma^{-1}(Y)=0$ and $\gamma'(0)$ is not tangent to $Y$. We have the homotopy-commutative square
\begin{equation}\label{manifoldpushout2}
\xymatrix{
P(Y,X)\ar[r]\ar[d] & Y\ar[d]\\
X-Y\ar[r] & X\\
}
\end{equation}
in which the top and left maps are evaluation at $0$ and at $1$ respectively. 

There are equivalences 
\begin{equation}\label{equiv1}
\hocofiber(P(Y,X)\rightarrow Y)\rightarrow \hocofiber(S(Y;\nu)\rightarrow Y)=Y^\nu
\end{equation}
and
\begin{equation}\label{equiv2}
\hocofiber(P(Y,X)\rightarrow Y)\rightarrow \hocofiber(X-Y\rightarrow X).
\end{equation}
The logic is as follows: 

For (\ref{equiv1}) we use the map $P(Y,X)\rightarrow S(Y;\nu)$ that sends $\gamma$ to the projection of $\gamma'(0)$ in the direction perpendicular to $Y$, normalized to have unit length. It is  a map over $Y$ between two spaces fibered over $Y$, and it is an equivalence because for each point in $Y$ the map of fibers is an equivalence.

For (\ref{equiv2}) we need to see that the homotopy-commutative square (\ref{manifoldpushout2}) is a homotopy pushout, in the sense that the associated map from the homotopy colimit of 
$$X-Y\leftarrow P(Y,X)\rightarrow Y$$
to $X$ is an equivalence. After choosing a tubular neighborhood of $Y$ in $X$, one can map $S(Y;\nu)$ to $P(Y,X)$ by using radial paths perpendicular to $Y$. This map is an equivalence because it is a one-sided inverse to an equivalence. It follows that in showing that the square is a homotopy pushout we may consider instead the square

$$\xymatrix{
S(Y;\nu)\ar[r]\ar[d] & Y\ar[d]\\
X-Y\ar[r] & X\\
}
$$
But this comes down to considering the same strictly commutative square (\ref{manifoldpushout}) that we began with.

Note that although a tubular neighborhood was used in proving (\ref{equiv2}) to be an equivalence, the definitions of (\ref{equiv1}) and (\ref{equiv2}) did not use it. This is the point of introducing $P(Y,X)$.

Now for the function spaces:   
Again we will obtain equivalences 
$$
\hocofiber(P(\mathcal B-\mathcal V,\mathcal M-\mathcal V)\rightarrow \mathcal M-\mathcal V)
\rightarrow (\mathcal B-\mathcal V)^\nu
$$
(where $\nu$ now means the bundle $TN/(TP\oplus TQ)$ on $\mathcal B-\mathcal V$) and 
$$
\hocofiber(P(\mathcal B-\mathcal V,\mathcal M-\mathcal V)\rightarrow \mathcal M-\mathcal V)
\rightarrow \hocofiber(\mathcal M-\mathcal B\to \mathcal M-\mathcal V)
$$

We define the space $P(\mathcal B-\mathcal V,\mathcal M-\mathcal V)$. A point in it is a map $\gamma:[0,1]\rightarrow \mathcal M$ meeting the following conditions. Write $\gamma(t)=(f_t,g_t)$. The conditions are:
\begin{enumerate}
\item $\gamma$ is smooth in the sense that the adjoint maps $(t,x)\mapsto f_t(x)$ and $(t,x)\mapsto g_t(x)$ from $[0,1]\times P$ and $[0,1]\times Q$ to $N$ are smooth. 
\item For every $t>0$, $\gamma_t$ is in $\mathcal M-\mathcal B$, that is, $f_t(P)\cap g_t(Q)=\emptyset$.
\item $\gamma_0\in \mathcal B-\mathcal V$, that is, (3a) there is exactly one point $(x_0,z_0,y_0)\in P\times N\times Q$ such that $f_0(x)=z_0=g_0(y)$ and (3b) $df_0\oplus dg_0:T_{x_0}P\oplus T_{y_0}Q\rightarrow T_{z_0}N$ is injective.
\item $\gamma'(0)$ is not tangent to $\mathcal B-\mathcal V$, that is, the vector $f'_0(x_0)-g'_0(y_0)\in T_{z_0}(N)$ does not belong to the subspace $(d_{x_0}f)(T_{x_0}P)\oplus (d_{y_0}g)(T_{y_0}Q)$. Here $f'$ and $g'$ are derivatives with respect to $t$.
\end{enumerate}

Consider the homotopy-commutative square
$$\xymatrix{
P(\mathcal B-\mathcal V,\mathcal M-\mathcal V)\ar[r]\ar[d] & \mathcal B-\mathcal V\ar[d]\\
\mathcal M-\mathcal B\ar[r] & \mathcal M-\mathcal V\, ,\\
}
$$
where the upper map and the left map take $\gamma=(f,g)$ to $(f_0,g_0)$ and $(f_1,g_1)$ respectively.
We argue much as in the finite-dimensional case. 

First, there is an equivalence $P(\mathcal B-\mathcal V,\mathcal M-\mathcal V)\rightarrow S(\mathcal B-\mathcal V;\nu)$ that respects the projection to $ \mathcal B-\mathcal V$, namely the map that takes $\gamma=(f,g)$ to the unit vector in $T_{z_0}N/(T_{x_0}P\oplus T_{y_0}Q)$ determined by the element $f'_0(x_0)-g'_0(y_0)$ of $T_{x_0}P\oplus T_{y_0}Q$. It is an equivalence because it is a map between spaces fibered over $\mathcal B-\mathcal V$ and it induces equivalences fiber by fiber.

Second, the square is a homotopy pushout. For this step, instead of trying to come up with a tubular neighborhood we reduce to the finite-dimensional case.

To show that the map from the homotopy colimit of 
$$
\mathcal M-\mathcal B\leftarrow P(\mathcal B-\mathcal V,\mathcal M-\mathcal V)\rightarrow \mathcal B-\mathcal V
$$
to $\mathcal M-\mathcal V$ is surjective on homotopy groups, let $X=S^k$ and take any map $\phi:X\rightarrow \mathcal M-\mathcal V$, with adjoint $\Phi=(F,G)$, $F:X\times P\rightarrow N$, $G:X\times Q\rightarrow N$. Deforming by a homotopy that stays within $\mathcal M-\mathcal V$, make $\Phi$ \lq\lq transverse to $\mathcal B-\mathcal V$\rq\rq\  in the sense that $F$ and $G$ together give a map $X\times P\times Q\rightarrow N\times N$ which is transverse to the diagonal. The preimage of the diagonal in $X\times P\times Q$ is a submanifold, and it is embedded in $X$ by the projection. Call its image $Y$. The normal bundle of $Y$ in $X$ is the pullback of $TN/(TP\oplus TQ)$ by $\phi$.

Now inverting the equivalence 
$$
\hocolim(X-Y\leftarrow P(Y,X)\rightarrow Y)\rightarrow X
$$
and composing with the obvious map
$$
\hocolim(X-Y\leftarrow P(Y,X)\rightarrow Y)\rightarrow 
\hocolim(\mathcal M-\mathcal B\leftarrow P(\mathcal B-\mathcal V,\mathcal M-\mathcal V)\rightarrow \mathcal M-\mathcal V)
$$
we get 
$$ 
X\rightarrow \hocolim(\mathcal M-\mathcal B\leftarrow P(\mathcal B-\mathcal V,\mathcal M-\mathcal V)\rightarrow \mathcal M-\mathcal V)
$$
a  lifting (up to homotopy) of $\phi$.
Essentially the same argument serves to lift a homotopy and prove the injectivity.

Taking homotopy fibers over $\mathcal M-\mathcal V$ all around, we obtain the needed equivalence F.

\section{Acknowledgments}

The authors would like to thank Harvard University for their hospitality.

\bibliographystyle{amsplain}

\bibliography{Bibliography}

\end{document}